\pdfoutput=1
\documentclass[10pt]{article}
\usepackage[a4paper, margin=1in]{geometry}
\usepackage{changepage}
\usepackage{amsmath, amssymb, amsthm}
\usepackage[utf8]{inputenc}
\usepackage[english]{babel}
\usepackage[numbers]{natbib}
\usepackage{url}
\usepackage[usenames]{color}
\usepackage{mathabx}
\usepackage[colorlinks=true]{hyperref}
\usepackage[capitalise]{cleveref}
\usepackage{tabularx}
\usepackage{breqn}
\usepackage[title]{appendix}
\usepackage{tikz}
\usepackage{xcolor}
\usepackage{comment}

\newtheoremstyle{mystyle}%                % Name
  {}%                                     % Space above
  {}%                                     % Space below
  {\itshape}%                             % Body font
  {}%                                     % Indent amount
  {\bfseries}%                            % Theorem head font
  {.}%                                    % Punctuation after theorem head
  { }%                                    % Space after theorem head, ' ', or \newline
  {\thmname{#1}\thmnumber{ #2}\thmnote{ (#3)}}%                                     % Theorem head spec (can be left empty, meaning `normal')
\theoremstyle{mystyle}
\newtheorem{theorem}{Theorem}

\newtheorem{lemma}[theorem]{Lemma}

\newtheorem{cor}[theorem]{Corollary}

\newtheorem{problem}[theorem]{Problem}

\newcommand{\floor}[1]{\left\lfloor #1 \right\rfloor}

\title{A Minimal Substitution Basis for the Kalmár Elementary Functions}

\author{Mihai Prunescu \footnote{Research Center for Logic, Optimization and Security (LOS) (Faculty of Mathematics and Computer Science, University of Bucharest, Academiei 14, Bucharest (RO-010014), Romania), e-mail address: {\tt mihai.prunescu@gmail.com}.} \footnote{Simion Stoilow Institute of Mathematics of the Romanian Academy (Research unit 5, P. O. Box 1-764, Bucharest (RO-014700), Romania), e-mail address: {\tt mihai.prunescu@imar.ro}.}, Lorenzo Sauras-Altuzarra \footnote{Kurt Gödel Society, Vienna, Austria, e-mail address: {\tt lorenzo@logic.at}.} \footnote{Institute for Logic and Data Science, Bucharest, Romania.}, Joseph M. Shunia \footnote{Independent researcher, Ann Arbor, Michigan, United States, e-mail address: {\tt jshunia@gmail.com}.}}

\date{May 2025 \\ \footnotesize{Revised: January 2026}}

\begin{document}

\maketitle

\begin{abstract} \noindent
We show that the class of Kalmár elementary functions can be inductively generated from the addition, the integer remainder, and the base-two exponentiation, hence improving previous results by Marchenkov and Mazzanti. We also prove that the substitution basis defined by these three operations is minimal. Furthermore, we discuss alternative substitution bases under arity constraints.
\\[2 mm]
\\[2 mm]
\\[2 mm]
\textbf{Note:} this is a preprint, the final version will appear in the Journal of Logic and Computation.
\end{abstract}

\section{Introduction} \label{introduction}

In this work, we denote by $ \mathbb{N} $ the set of non-negative integers, and by $ \dotdiv $ the \textbf{truncated subtraction}: $ x \dotdiv y := \max( x - y , 0 ) $. In addition, given an integer $ r \geq 1 $ and a $ r $-ary operation $ f $ on $ \mathbb{N} $, we will often write $ f $ as $ ( x_1 , \ldots , x_r ) \mapsto f ( x_1 , \ldots , x_r ) $ or even just as $ f ( x_1 , \ldots , x_r ) $. For example, we denote the addition as $ x + y $. We also apply the conventions $ 0^0 := 1 $ and $ \floor{x / 0} := 0 $, so, in particular, by applying the definition of the \textbf{integer remainder}, $ x \bmod y := x \dotdiv y \floor{x / y} $, we have that $ x \bmod 1 = 0 $ and $ x \bmod 0 = x $.

Given a set $ \mathcal{B} $ of operations on $ \mathbb{N} $, we denote by $ \langle \mathcal{B} \rangle $ the closure under substitution of the set $ \mathcal{B} \cup C \cup P $, where $ C $ is the set of constant operations on $ \mathbb{N} $ (which can be identified with non-negative integers) and $ P $ is the set of projections on $ \mathbb{N} $ (which can be identified with variables). And we say that $ \mathcal{B} $ is a \textbf{substitution basis} for a given set $ G $ of operations on $ \mathbb{N} $ if, and only if, $ \langle \mathcal{B} \rangle = G $ (cf. Mazzanti \cite[Section 2.2]{mazzanti2002plainbases}). For instance, if $ S := \langle x + y , \ x^y \rangle $, then $ x^{z + 2} \in S $: indeed, the projection $ ( x , y , z ) \mapsto z $ and the constant operation $ ( x , y , z ) \mapsto 2 $ belong to $ S $, so, because of the closure under substitution of $ S $, the operation $ ( x , y , z ) \mapsto x^{z + 2} $ also belongs to $ S $. Observe that some authors do not admit the constant operations as basic elements in a substitution basis (cf. Marchenkov \cite[Introduction]{marchenkov2007superposition}), but, as we are already admitting any variable $ x_k $, it is not unreasonable to also admit any constant $ k $.

The \textbf{Kalmár elementary functions} are the elements of the smallest set of operations on $ \mathbb{N} $; containing the constants, the projections, the addition, and the truncated difference; that is closed under substitution, bounded summation, and bounded multiplication (cf. Clote \cite[Definition 3.64]{Clote1999} and Cutland \cite[Definition 4.1]{cutland1980}). They are often characterized as the operations on $ \mathbb{N} $ that can be computed in iterated exponential time (cf. Cutland \cite[Theorem 4.11]{cutland1980}, Marchenkov \cite[Introduction]{marchenkov2007superposition}, and Oitavem \cite[Introduction]{oitavem1997}); and they form the class $ \mathcal{E}^3 $ in the Grzegorczyk hierarchy (see Grzegorczyk \cite{grzegorczyk1953someclasses}), which we will denote just by $ \mathcal{E} $, and which encompasses most of the operations on $ \mathbb{N} $ that frequently arise in mathematical contexts (cf. Campagnolo et al. \cite[Introduction]{CampagnoloEtAl}). In practical terms, a significant portion of the real-world computation can be faithfully modeled within the class of Kalmár elementary functions, since the corresponding code avoids unbounded recursion.

In 1964, Rödding \cite{Rodding1964} (in German) showed that a finite quantity of elementary operations suffices to form a substitution basis for the class of Kalmár elementary functions, which marked the beginning of the search for substitution bases consisting of ``simple'' arithmetic operations. This search was considered concluded in 2002 (see Marchenkov \cite[Introduction]{marchenkov2007superposition}), when Mazzanti \cite{mazzanti2002plainbases} proved that \begin{equation} \label{MazzantiEq} \langle x + y , \ x \dotdiv y , \ x y , \ \floor{x / y} , \ x^y \rangle = \mathcal{E} . \end{equation}

For reference, we mention that recently we have obtained explicit representations of some classical number-theoretic functions using the substitution basis from Mazzanti's \cref{MazzantiEq}; such as the factorial function $ n ! $ (see Prunescu \& Sauras-Altuzarra \cite{prunescusauras2024factorial}), the sums of binomial coefficients (see Shunia \& Sauras-Altuzarra \cite[Theorem 3.1]{shuniasauras2025}), the $n$-th prime function $ p_n $ (see Prunescu \& Shunia \cite[Theorem 10.1]{prunescushunia2024primes}), Euler's totient function $\varphi(n)$ (see Prunescu \& Sauras-Altuzarra \cite[Theorem 4.12]{prunescusauras2025manyfunctions}), the greatest common divisor $ \gcd ( m , n ) $ (see Prunescu \& Shunia \cite{prunescushunia2024gcd}), and any C-recursive integer sequence with rational coefficients (see Prunescu \cite{prunescu2025fibonacci2}, Prunescu \& Sauras-Altuzarra \cite{prunescusauras2025fibonacci}, and Prunescu \& Shunia \cite{prunescushunia2025fibonacci3}).

Slightly after Mazzanti's \cref{MazzantiEq}, in 2007, Marchenkov \cite[Corollary 2]{marchenkov2007superposition} established that \begin{equation} \label{MarchenkovEq} \langle x + y , \ x \bmod y , \ x^2 , \ 2^x \rangle = \mathcal{E} . \end{equation}

The present article improves Marchenkov's result by eliminating the squaring operation from the above substitution basis. Concretely, \cref{MainTheorem} asserts that $$ \langle x + y , \ x \bmod y , \ 2^x \rangle = \mathcal{E} ; $$ and \cref{SectionMinimality} is devoted to show that this substitution basis is minimal: if $ F \subsetneqq \{ x + y , x \bmod y , 2^x \} $, then $ \langle F \rangle \subsetneqq \mathcal{E} $. In addition, \cref{SectionBadSets} discusses why certain sets cannot serve as substitution bases for the class of Kalmár elementary functions. And, finally, \cref{SectionLowArity} presents one finite substitution basis for the class of unary Kalmár elementary functions that consists of unary operations only, and two substitution bases for the whole class of Kalmár elementary functions that consist of exactly one binary operation.

\section{A substitution basis for the Kalmár elementary functions}

We start by showing that any Kalmár elementary function can be represented in terms of the addition, the integer remainder, and the base-two exponentiation.

\begin{theorem} \label{MainTheorem} $ \langle x + y , \ x \bmod y , \ 2^x \rangle = \mathcal{E} $. \end{theorem}

\begin{proof} By Marchenkov's \cref{MarchenkovEq}, it is enough to prove that $ x^2 \in \langle x+y, x \bmod y, 2^x \rangle $.

Since $0 \leq x^2 < 2^x + x$ and $2^{x + x} - k(2^x + x) = x^2$, where $k = 2^x - x$, then $2^{x+x} \bmod (2^x + x) = x^2$. \end{proof}

Now, we provide a proof of \cref{thmimportantoperations}, which asserts that the truncated subtraction, the product, and the integer division certainly can be represented by using the substitution basis from Marchenkov's \cref{MarchenkovEq}. This result relies in \cref{lemmadifference} and \cref{lemmaquotient}, identities applied by Marchenkov \cite{marchenkov2007superposition} but whose proof appeared to be missing.

\begin{lemma} \label{lemmadifference}
If $ x $ and $ y $ are non-negative integers, then $$ x \dotdiv y = \left( (2^{x+y} + x) \bmod (2^{x+y} + y) \right) \bmod (2^{x+y}+x) . $$
\end{lemma}

\begin{proof} It is clear that, if $ x \geq y $, then $ 0 \leq x - y < 2^{x+y} + y $; so the division with remainder reads $ 2^{x+y} + x = 1 \cdot (2^{x+y} + y) + (x-y) $.

Also, in the case that $ x < y $, we have that $ 0 \leq 2^{x+y} + x < 2^{x+y} + y $; so the division with remainder has the form $ 2^{x+y} + x = 0 \cdot (2^{x+y} + y) + (2^{x+y} + x) $.

Altogether $$ (2^{x+y}+x) \bmod (2^{x+y}+y) = \begin{cases} x-y & \textup{ if } x \geq y , \\ 2^{x+y}+x & \textup{ if } x < y \end{cases} $$ and, by applying that $ 0 \leq x-y < 2^{x+y}+x $ if $ x \geq y $, the statement follows. \end{proof}

\begin{lemma} \label{lemmaquotient} If $ x $ and $ y $ are non-negative integers, then $$ \floor{x / y} = \left( 2 ( x + 1 ) \left( x \dotdiv ( x \bmod y ) \right) \right) \bmod ( 2 (x+1) y \dotdiv 1 ) . $$ \end{lemma}

\begin{proof} If $ y = 0 $, then, by the convention mentioned in \cref{introduction}, we have that $ \floor{x / y} = \floor{x / 0} = 0 = 0 \bmod 0 = \left( 2 ( x + 1 ) \left( x \dotdiv x \right) \right) \bmod ( 0 \dotdiv 1 ) = \left( 2 ( x + 1 ) \left( x \dotdiv ( x \bmod 0 ) \right) \right) \bmod ( 2 (x+1) 0 \dotdiv 1 ) = \left( 2 ( x + 1 ) \left( x \dotdiv ( x \bmod y ) \right) \right) \bmod ( 2 (x+1) y \dotdiv 1 ) $. Otherwise, by applying that $ x \dotdiv ( x \bmod y ) = y \lfloor x / y  \rfloor $, we have that $ 2 ( x + 1 ) \left( x \dotdiv ( x \bmod y ) \right) = 2 x y \lfloor x / y \rfloor + 2 y \lfloor x / y \rfloor = \lfloor x / y \rfloor ( 2 x y + 2 y \dotdiv 1 ) + \lfloor x / y \rfloor $ and thus, by applying that $ 0 \leq \lfloor x / y \rfloor \leq x < 2 x y + 2 y \dotdiv 1 = 2 ( x + 1 ) y \dotdiv 1 $, the statement follows. \end{proof}

\begin{theorem} \label{thmimportantoperations} The operations $ x \dotdiv y $, $ x y $, and $ \floor{x / y} $ belong to $ \langle x + y , \ x \bmod y , \ 2^x \rangle $. \end{theorem}

\begin{proof} Let $ S = \langle x + y , \ x \bmod y , \ 2^x \rangle $. \cref{lemmaquotient} implies that $ \floor{x / y} \in \langle x + y , \ x \dotdiv y , \ x \bmod y , \ 2 x y \rangle $. And, by applying \cref{lemmadifference}, we get $ x \dotdiv y \in S $. Thus, by applying the fact that $x^2 \in S$ (see the proof of \cref{MainTheorem}), we have that $ 2 x y = (x+y)^2 \dotdiv (x^2 + y^2) \in S $ and hence $ \floor{x / y} \in S $. Finally, notice that $ x y = \floor{2 x y / 2} \in S $. \end{proof}

\section{Minimality of the main substitution basis} \label{SectionMinimality}

In \cref{MainTheorem}, we have seen that the set $ \{ x + y, x \bmod y, 2^x \} $ is a substitution basis for the Kalmár elementary functions. In this section we prove that it is, in addition, \textit{minimal}: $ x + y \notin \langle x \bmod y , \ 2^x \rangle $, $ x \bmod y \notin \langle x + y , \ 2^x \rangle $, and $ 2^x \notin \langle x + y , \ x \bmod y \rangle $.

Recall that the definition of substitution basis that we gave in \cref{introduction} allows the usage of constants when generating functions from a given basis.

\begin{lemma} \label{lemmamodexp} If $ t ( x ) \in \langle x \bmod y , \ 2^x \rangle $, then there is a non-negative integer $ B $ such that, for every non-negative integer $ a $, the number $ t ( a ) $ is a power of two or $ t ( a ) \leq \max ( B , a ) $. \end{lemma}

\begin{proof} The proof goes by induction on the elements of $ \langle  x \bmod y , \ 2^x \rangle $.

The statement is obvious for the constants and the variable $ x $.

Now, let $ t_1 ( x ) $ and $ t_2 ( x ) $ be elements of $ \langle x \bmod y , \ 2^x \rangle $ that satisfy the property given in the statement: there are two non-negative integers $ B_1 $ and $ B_2 $ such that, for every non-negative integer $ a $, the number $ t_1 ( a ) $ is a power of two or $ t_1 ( a ) \leq \max ( B_1 , a ) $, while the number $ t_2 ( a ) $ is a power of two or $ t_2 ( a ) \leq \max ( B_2 , a ) $.

If $ t ( x ) = 2^{t_1 ( x )} $, then, for every non-negative integer $ a $, the number $ t ( a ) $ is a power of two, and it does not matter how we choose $ B $ in this case.

If $ t ( x ) = t_1 ( x ) \bmod t_2 ( x ) $, then we can set $ B = \max ( B_1 , B_2 ) $: indeed, for every non-negative integer $ a $, we have that:
\begin{enumerate}
    \item if $ t_2 ( a ) = 0 $, then $ t ( a ) = t_1 (a ) \bmod 0 = t_1 ( a ) $ and consequently $ t ( a ) $ is a power of two or $ t ( a ) \leq \max ( B_1 , a ) \leq \max ( B , a ) $;
    \item if $ 0 < t_2 ( a ) \leq \max ( B_2 , a ) $, then $ t ( a ) < t_2 ( a ) \leq \max ( B_2 , a ) \leq \max ( B , a ) $;
    \item if $ t_1 ( a ) \leq \max( B_1 , a ) $, then $ t ( a ) \leq t_1 ( a ) \leq \max ( B_1 , a ) \leq \max ( B , a ) $; and
    \item if $ t_1 ( a ) $ and $ t_2 ( a ) $ are powers of two, then $ t ( a ) $ is a power of two or $ t ( a ) = 0 \leq \max ( B , a ) $. \qedhere
\end{enumerate} \end{proof}

\begin{theorem}\label{theononaddition} $ x + y \notin \langle x \bmod y , \ 2^x \rangle $. \end{theorem}

\begin{proof} Suppose the contrary. It is then clear that $ x + 1 \in \langle x \bmod y , \ 2^x \rangle $. According to \cref{lemmamodexp}, there is some non-negative integer $ B $ such that, for every non-negative integer $ a $, the number $ a + 1 $ is a power of two or $ a + 1 \leq \max ( B , a ) $. Therefore, every integer $ a > B $ should be a power of two, which is false. \end{proof}

\begin{lemma} \label{lemmastrictlyincreasing} If $ t ( x ) \in \langle x + y , \ 2^x \rangle $ is non-constant, then $ t ( x ) $ is strictly increasing. \end{lemma}

\begin{proof} The proof goes by induction on the elements of $ \langle x + y , \ 2^x \rangle $. The statement is obvious for the variable $ x $. And, if $ t_1 ( x ) $ and $ t_2 ( x ) $ are strictly increasing elements of $ \langle x + y , \ 2^x \rangle $, then $ t_1 (x) + t_2 (x) $ and $ 2^{t_1 ( x )} $ are also strictly increasing. \end{proof}

\begin{theorem} $ x \bmod y \notin \langle x + y , \ 2^x \rangle $. \end{theorem}

\begin{proof} Suppose the contrary. Then the operation $ x \bmod 2 $, which is not strictly increasing, is also a non-constant element of $ \langle x + y , \ 2^x \rangle $, in contradiction with \cref{lemmastrictlyincreasing}. \end{proof}

\begin{lemma} \label{lemmalineargrowth} If $ t ( x ) \in \langle x + y , \ x \bmod y \rangle $, then $ t ( x ) < A x + B $ for some non-negative integers $ A $ and $ B $. \end{lemma}

\begin{proof} The proof goes by induction on the elements of $ \langle x + y , \ x \bmod y \rangle $. The statement is obvious for the constants and the variable $ x $. And, if $ t_1 ( x ) $ and $ t_2 ( x ) $ are elements of $ \langle x + y , \ x \bmod y \rangle $ such that $ t_1 ( x ) < A_1 x + B_1 $ and $ t_2 ( x ) < A_2 x + B_2 $ for some non-negative integers $ A_1 $, $ A_2 $, $ B_1 $, and $ B_2 $, then it is clear that $ t_1 (x) + t_2 (x) < ( A_1 + A_2 ) x + ( B_1 + B_2 ) $ and $ t_1 ( x ) \bmod t_2 ( x ) \leq t_1 ( x ) < A_1 x + B_1 $. \end{proof}

\begin{theorem} $ 2^x \notin \langle x + y , \ x \bmod y \rangle $. \end{theorem}

\begin{proof} Suppose the contrary. Then, \cref{lemmalineargrowth} yields the existence of two non-negative integers $ A $ and $ B $ such that $ 2^x < A x + B $, which is clearly impossible. \end{proof}

\section{Some sets which are not substitution bases} \label{SectionBadSets}

In this section we show that, if in the basis $\{x + y, x \bmod y, 2^x\}$, the operation $x+y$ is replaced with the operation $2x$, $x+1$, or $x \dotdiv y$, then the resulting set is not a substitution basis for the class of Kalmár elementary functions. We start with \cref{lemmatwomodexp}, which is similar with \cref{lemmamodexp}.

\begin{lemma} \label{lemmatwomodexp} If $ t ( x ) \in \langle 2x, \, x \bmod y , \ 2^x \rangle $, then there are two non-negative integers $ A $ and $ B $ such that, for every non-negative integer $ a $, the number $ t ( a ) $ is a power of two or $ t ( a ) \leq A \max ( B , a ) $. \end{lemma}

\begin{proof} The proof goes by induction on the elements of $ \langle 2x,\, x \bmod y , \ 2^x \rangle $.

The statement is obvious for the constants and the variable $ x $.

Now, let $ t_1 ( x ) $ and $ t_2 ( x ) $ be elements of $ \langle 2x, \, x \bmod y , \ 2^x \rangle $ that satisfy the property given in the statement: there are four non-negative integers $ A_1 $, $ A_2 $, $ B_1 $, and $ B_2 $ such that, for every non-negative integer $ a $, the number $ t_1 ( a ) $ is a power of two or $ t_1 ( a ) \leq A_1 \max( B_1 , a )$, while the number $ t_2 ( a ) $ is a power of two or $ t_2 ( a ) \leq A_2 \max( B_2 , a ) $.

If $ t ( x ) = 2 t_1 ( x ) $, then, for every non-negative integer $ a $, the number $ t ( a ) $ is a power of two or $ t ( a ) \leq 2 A_1 \max ( B_1 , a ) $, so we can set $ A = 2 A_1 $ and $ B = B_1 $.

If $ t ( x ) = 2^{t_1 ( x )} $, then, for every non-negative integer $ a $, the number $ t ( a ) $ is a power of two, and it does not matter how we choose $ A $ and $ B $ in this case.

If $ t ( x ) = t_1 ( x ) \bmod t_2 ( x ) $, then we can take $ A = \max ( A_1 , A_2 ) $ and $ B = \max ( B_1 , B_2 ) $: indeed, for every non-negative integer $ a $, we have that:
\begin{enumerate}
    \item if $ t_2 ( a ) = 0 $, then $ t ( a ) = t_1 ( a ) \bmod 0 = t_1 ( a ) $ and consequently $ t ( a ) $ is a power of two or $ t ( a ) \leq A_1 \max ( B_1 , a ) \leq A \max ( B , a ) $;
    \item if $ 0 < t_2 ( a ) \leq A_2 \max( B_2 , a ) $, then $ t ( a ) < t_2 ( a ) \leq A_2 \max( B_2 , a ) \leq A \max ( B , a ) $;
    \item if $ t_1 ( a ) \leq A_1 \max ( B_1 , a ) $, then $ t ( a ) \leq t_1 ( a ) \leq A_1 \max ( B_1 , a ) \leq A \max ( B , a ) $; and
    \item if $ t_1 ( a ) $ and $ t_2 ( a ) $ are powers of two, then $ t ( a ) $ is a power of two or $ t ( a ) = 0 \leq A \max ( B , a ) $. \qedhere
\end{enumerate} \end{proof}

\begin{theorem} \label{thmweakset} $ \langle 2 x , \ x \bmod y , \ 2^x \rangle \subsetneqq \mathcal{E} $. \end{theorem}

\begin{proof} Consider the function $ x^2 $, which belongs to $ \mathcal{E} $, and suppose that $x^2 \in \langle 2 x , \ x \bmod y , \ 2^x \rangle $. According to \cref{lemmatwomodexp}, there are two non-negative integers $ A $ and $ B $ such that, for every non-negative integer $ a $, the number $ a^2 $ is a power of two or $ a^2 \leq A \max ( B , a ) $. Therefore, for every integer $ a > \max ( A , B ) $, $ a^2 $ should be a power of two, which is false. \end{proof}

In a previous version of this article, \cref{theoremweakset2} was stated as a question, but Emil Jeřábek provided a proof (pers. comm.).

\begin{theorem} \label{theoremweakset2} $ \langle x + 1, \ x \bmod y , \ 2^x \rangle \subsetneqq \mathcal{E} $. \end{theorem}

\begin{proof} Let $ B $ denote the set of Kalmár elementary functions $ f ( \vec{x} ) $ such that, if $ r $ is the arity of $ f ( \vec{x} ) $, then there exists a computable function $ g ( x ) $ and a number $ i \in \{ 1 , \ldots , r \} $  such that $ f ( \vec{x} ) \leq g ( x_i ) $. We easily observe that the operations $ x + 1 $, $ x \bmod y $, and $ 2^x $ belong to $ B $, so, as clearly $ \langle B \rangle = B $, we get $ \langle x \bmod y , \ x + 1 , \ 2^x \rangle \subseteq B $. However, it is straightforward to see that $ x + y \in \mathcal{E} \backslash B $. \end{proof}

\begin{cor} $ \langle x \dotdiv y, \ x \bmod y , \ 2^x \rangle \subsetneqq \mathcal{E} $. \end{cor}

\begin{proof} The proof is the same as that of \cref{theoremweakset2}, but replacing $ x + 1 $ with $ x \dotdiv y $ in the set. \end{proof}

Note that the proof of \cref{theoremweakset2} could have been also adapted in order to prove \cref{theononaddition} and \cref{thmweakset}.

Finally, we state a research question, \cref{ProblemAlternativeBasis}.

\begin{problem} \label{ProblemAlternativeBasis} Is $ \langle x+y,\, \lfloor x / y \rfloor,\, 2^x \rangle $ equal to $ \mathcal{E} $? \end{problem}

Intuitively, as all three operations displayed in \cref{ProblemAlternativeBasis} are ultimately monotonic in both arguments, it seems to be challenging to generate periodic behavior (as required for $x \bmod y$) solely from these. But, on the other hand, there are cases in which it is possible: consider, for example, the identity
$$ 2^{x \bmod 2} = \left\lfloor \frac{2^x}{2^{\lfloor x/2 \rfloor + \lfloor x/2 \rfloor}} \right\rfloor , $$ which is not hard to prove.

\section{Other finite substitution bases with low-arity operations} \label{SectionLowArity}

We conclude the article by discussing the existence of substitution bases for the class of Kalmár elementary functions that consist of only two operations, or even of just one operation.

First of all, we observe that it is not difficult to construct \textit{artificial} substitution bases consisting in just one ternary operation that is not really mathematically relevant, for instance that defined as $$ \begin{cases} 2^x, & \textrm{ if } z = 0 , \\ x+y, & \textrm{ if } z = 1 , \\ x \bmod y, & \textrm{ if } z \geq 2 . \end{cases} $$ Therefore, we should restrict our attention to unary and binary operations only.

Secondly, we note that a finite set of unary operations cannot be a substitution basis for the class of Kalmár elementary functions, because of the simple reason that a term which is built up from unary function-symbols can contain at most one free variable. However, the question of whether it can generate the class of \textit{unary} Kalmár elementary functions is more reasonable.

\textbf{Cantor’s pairing function}, which we denote by $ [ x \mid y ] $, is a bijection that maps each pair $ ( x , y ) $ of non-negative integers into the non-negative integer $ ( x + y ) ( x + y + 1 ) / 2 + x $. In addition, we consider in this section the unary operations $ l ( x ) $ and $ r ( x ) $ such that $ l ( [ x \mid y ] ) = x $ and $ r ( [ x \mid y ] ) = y $. Note that $[x \mid y]$, $l(z)$ and $r(z)$ are Kalmár elementary functions; and $$ [l(x) \mid r(x)] = x $$ for every $x \in \mathbb N$.

\begin{theorem} \label{MarchenkovThm} Let $ \Phi $ be a class of operations in $ \mathbb{N} $, containing $ [ x \mid y ] $, $ l ( x ) $, and $ r ( x ) $; and closed under substitution. In addition, suppose that $ \Phi $ has a finite substitution basis. Then the class of unary operations in $ \Phi $ has a substitution basis consisting of two unary operations only. If, moreover, $ \Phi $ contains the operations $ x + 1 $ and $ x \dotdiv y $, then $ \Phi $ has a substitution basis consisting of one binary operation only. \end{theorem}

Notice that \cref{MarchenkovThm}, which is due to Marchenkov \cite[Theorem 2]{marchenkov1991} (in Russian), answers affirmatively the questions posed above, as we can clearly instantiate $ \Phi $ to $ \mathcal{E} $. We now present alternative proofs of parts of this result.

\begin{theorem} \label{theo:unarybasis} There is a finite substitution basis for the class of unary Kalmár elementary functions consisting of unary operations only. \end{theorem}

\begin{proof}[Proof (due to Emil Jeřábek, pers. comm.).] Let $ \mathcal{B} $ be a finite substitution basis for $\cal E$ consisting of unary and binary operations only. In particular, we can suppose that the chosen basis generates all constants (which is true, for instance, for the basis from \cref{MainTheorem}, because $ 0 = x \bmod x $, $ 1 = 2^0 $, $ 2 = 1 + 1 $, etc.).

The functions $$d(x) = [x \mid x], \,\,\,\, u(x) = [l(x) \mid x], \,\,\,\, v(x) = [l(r(x)) \mid [l(x) \mid r(r(x))]]$$ are in ${\cal E} = \langle \mathcal{B} \rangle$, and given two functions $f,h \in \langle \mathcal{B} \rangle$ with $f$ unary and $h$ binary, let $$f'(x) = [f(l(x)) \mid r(x)], \,\,\,\, h''(x) = [h(l(x)), l(r(x)) ) \mid r(r(x))].$$ Observe that, if $f, g \in \langle \mathcal{B} \rangle$ are unary, then $(g \circ f)' = g' \circ f'$. The set $$\mathcal{U} = \{l, r, d, u, v\} \cup \{f' : f \in \mathcal{B} \textrm{ is a unary function}\} \cup \{h'' : h \in \mathcal{B} \textrm{ is a binary function} \}$$ is a collection of unary functions in ${\cal E}$, and so is $\langle \mathcal{U} \rangle$, the family of functions obtained from $ \mathcal{U} $, the constants, and the identity function. Since $F(x) = l(F'(d(x)))$, it follows that $F \in \langle \mathcal{U} \rangle \Leftrightarrow F' \in \langle \mathcal{U} \rangle$, for any unary function $F$.

We claim that $\mathcal{U}$ is a finite basis for the unary Kalmár elementary functions. It suffices to prove by induction on the complexity of a unary $F \in \langle \mathcal{B} \rangle$ that $F \in \langle \mathcal{U} \rangle$.

The base case is immediate: if $F$ is the identity, then it belongs to $\langle \mathcal{U} \rangle$.

Suppose $F = g \circ f$ with $f \in \langle \mathcal{B} \rangle $  and $g \in \mathcal{B}$, both unary. Then $F' = (g \circ f)' = g' \circ f'$ and $f \in \langle \mathcal{U} \rangle$ by the induction assumption, and hence $f' \in \langle \mathcal{U} \rangle$, and $g' \in \mathcal{U}$ by the definition of $\mathcal{U}$, so that $F' \in \langle \mathcal{U} \rangle$, and hence $F \in \langle \mathcal{U} \rangle$.

Suppose $F(x) =h(f(x), g(x)) $ with $f,g \in \langle \mathcal{B} \rangle$ unary and $h \in \mathcal{B}$ binary. By inductive assumption $f, g \in \langle \mathcal{U} \rangle$ and hence $f', g' \in \langle \mathcal{U} \rangle$, and using the definitions of $u$, $v$ and the definition of $f'$ and $g'$, we have that, for all $x$, $$(g' \circ u)(x) = [ g(l(x)) \mid x],$$ $$(v \circ g' \circ u) (x) = [l(x) \mid [g(l(x)) \mid r(x) ]],$$ $$(f' \circ v \circ g' \circ u)(x) = [f(l(x)) \mid [g(l(x)) \mid r(x) ]],$$ $$(h'' \circ f' \circ v \circ g' \circ u)(x) = [h (f ( l(x) ), g(l(x))) \mid r(x) ].$$ It is easy to check that $F'(x) = [h(f(l(x)), g(l(x))) \mid r(x)]$, so $F' \in \langle \mathcal{U} \rangle$, and therefore $F \in \langle \mathcal{U} \rangle$. \end{proof}

For example, let $ f ( x ) = 2 x $. Then we can represent the operation $ f ( x ) $ in the basis for the class of unary Kalmár elementary functions that is given in the proof of \cref{theo:unarybasis} by setting $ a(x) = [ l(x) + l(r(x)) \mid r(r(x)) ] $: indeed, by performing a few calculations, we obtain $ f ( x ) = l ( [ f ( x ) \mid x ] ) = l ( [ x + x \mid x ] ) = l ( a ( [ x \mid [ x \mid x ] ] ) ) = l ( a ( u ( [ x \mid x ] ) ) ) = l ( a ( u ( d ( x ) ) ) ) $.

\begin{theorem}\label{theo:onlyonebinary} There is a basis for the class of Kalmár elementary functions that consists of one binary operation only. \end{theorem}

\begin{proof}[First proof (due to Emil Jeřábek, pers. comm.).] By applying \cref{theo:unarybasis}, there is a positive integer $ k $ and a substitution basis $ B $ for the class of unary Kalmár elementary functions that consists of exactly $ k $ unary operations $ f_0 ( x ) $, $ \ldots $, $ f_{k - 1} ( x ) $ such that $ f_0(x) = 2(x + 1) $.

In addition, it is easy to prove that, if $ i $ is a non-negative integer, then the $ i $-th iteration of $ f_0 ( x ) $, denoted, as usual, as $ f_0^i ( x ) $, is equal to $ 2^i x + 2^{i + 1} \dotdiv 2 $.

Now, let $$ g (x, y) = \begin{cases} f_i(x) & \text{if } y = f_0^i(x) \text{ for some } i \in \{ 0 , \ldots , k - 1 \} , \\ [ x / 2 \mid ( y \dotdiv 1 ) / 2 ] & \text{if } x \text{ is even and } y \text{ is odd} , \\ 0 & \text{otherwise} , \end{cases} $$ which is well-defined because the cases are disjoint: indeed, notice that, if $ i $ is a positive integer, then every term of the sequence $ f_0^i ( x ) $ is even.

We have that $ f_0 ( x ) = g ( x , x ) $ and $ f_i ( x ) = g ( x , f_0^i ( x ) ) $ for every $ i \in \{ 1 , \ldots , k - 1 \} $, so the set $ \langle g ( x , y ) \rangle $ contains all unary Kalmár elementary functions.

Thus, $ 2 x $ and $ 2 x + 1 $ belong to $ \langle g ( x , y ) \rangle $, whence it follows that $ [ x \mid y ] = g ( 2 x , 2 y + 1 ) \in \langle g ( x , y ) \rangle $.

Finally, let $ s ( x ) = l ( x ) + r ( x ) $ and $ m(x) = l(x) \bmod r(x) $.

The operations $ s ( x ) $, $ m ( x ) $, and $ 2^x $ are unary Kalmár elementary functions, so they belong to $ \langle g(x,y) \rangle $ and consequently we have that $ x + y = s ( [ x \mid y ] ) \in \langle g ( x , y ) \rangle $ and $ x \bmod y = m([x \mid y]) \in \langle g(x,y) \rangle $. Hence, by applying \cref{MainTheorem}, we deduce that $ \langle g ( x , y ) \rangle = \mathcal{E} $. \end{proof}

\begin{proof}[Second proof] Consider the binary Kalmár elementary function $ h ( x , y ) $ defined as $$ \begin{cases}
    2^x, & \textrm{if } x=y,\\
    3^{x+1}, & \textrm{if } y=2^x,\\
    5^{y+1}, & \textrm{if } x = 2^y,\\
    a+b, & \textrm{if } (x, y) = (3^{a+1}, 5^{b+1}),\\
    a \bmod b, & \textrm{if } (x, y) = (5^{a+1}, 3^{b+1}),\\
    0, & \textrm{otherwise.}
\end{cases}$$

Then, by applying \cref{MainTheorem} and the identities \[ 2^x = h(x,x) , \] \[ x + y = h(3^{x+1}, 5^{y+1}) = h(h(x,2^x), h(2^y,y)) = h(h(x,h(x,x)), h(h(y,y),y)) , \textrm{ and} \] \[ x \bmod y = h(5^{x+1}, 3^{y+1}) = h ( h ( 2^x , x ) , h ( y , 2^y ) ) = h ( h ( h ( x , x ) , x ) , h ( y , h ( y , y ) ) ) ; \] we conclude that the set $ \{ h ( x , y ) \} $ is a substitution basis for the class of Kalmár elementary functions. \end{proof}

Note that the substitution basis $ \{ h ( x , y ) \} $ from the last proof allows to generate every constant from one variable: indeed, $ 0 = x \bmod x = h( h( h(x,x), x) , h(x, h(x,x))) $, $ 1 = h(0,0) $, $ 2 = h(1,1) $, $ 3 = 1 + 2 $, and so on.

\section*{Acknowledgements}

We thank the anonymous referees for their invaluable, expert, and thoroughly careful review; Emil Jeřábek for his significant contributions; and Eugenio Omodeo for his attentive reading and kind comments.

\end{document}